\documentclass[12pt]{article}
\usepackage{mathtools}
\usepackage[utf8]{inputenc}
\usepackage{amssymb,amsmath,amsfonts,eucal,mathrsfs,amsthm} 
\usepackage[colorinlistoftodos]{todonotes}
\usepackage{xcolor}
\usepackage{enumitem}
\newtheorem{theorem}{Theorem}
\newtheorem{proposition}[theorem]{Proposition}
\newtheorem{lemma}[theorem]{Lemma}

\theoremstyle{definition}
\newcommand{\R}{\mathbb{R}}
\newcommand{\Z}{\mathbb{Z}}
\newcommand{\Q}{\mathbb{Q}}
\newcommand{\Sf}{\mathbb{S}}

\newcommand{\Ric}{\mbox{Ric}}

\newcommand{\trace}{\mbox{tr\,}}

\def\<{{\langle}}
\def\>{{\rangle}}
\def\B{\mathcal{B}}
\def\CP{\mathord{\mathbb C}\mathord{P}}

\def\n{\nabla}

\def\a{\alpha}

\def\be{\begin{equation} }
\def\ee{\end{equation} }

\newcommand\blfootnote[1]{\begingroup
\renewcommand\thefootnote{}\footnote{#1}
\addtocounter{footnote}{-1}
\endgroup}
\begin{document}

\title{Topological and geometric rigidity of 
nonnegatively curved submanifolds}
\author{Theodoros Vlachos}
\date{}
\maketitle

\begin{abstract}
We investigate the topology and geometry of compact 
submanifolds in space forms of nonnegative curvature 
satisfying a lower bound on the sectional 
curvature that depends only on the length of the mean 
curvature vector of the immersion. We show that this 
condition imposes strong constraints on either the 
topology or geometry of the submanifold. More precisely, 
we prove that, for dimensions $n\ge5$, the manifold 
is homeomorphic to a sphere, whereas in dimension four 
it is either diffeomorphic to $\Sf^4$ or isometric to 
$\CP^2$, in which case the immersion is the standard 
embedding into $\Sf^7$ followed by an umbilical inclusion. 
Additionally, we provide examples showing that the 
spherical conclusion cannot in general be 
strengthened to an extrinsic rigidity statement, 
since there exist many geometrically distinct spheres 
satisfying the pinching condition.
\end{abstract}
\blfootnote{\textup{2020} \textit{Mathematics Subject 
Classification}: 53C20, 53C40, 53C42.}
\blfootnote{\textit{Key words}:
Compact submanifolds, sectional, isotropic and 
mean curvature, homology groups, self-dual and 
anti-self-dual harmonic forms, 
Bochner-Weitzenb\"ock operator.}

\section{Introduction}
This paper focuses on the topological and 
geometric rigidity of compact submanifolds in space 
forms under a sectional curvature pinching condition. 
Previous works have studied compact submanifolds 
in space forms with lower bounds on their sectional 
curvature, establishing rigidity results under conditions 
such as minimality or parallel mean curvature vector 
field (see, for example, \cite{guxu, Itoh1, Itoh2, yau}).

In the absence of such geometric constraints, topological 
rigidity is the only result that can generally be anticipated 
under sectional curvature pinching conditions. Notably, 
significant results have been achieved by various 
researchers when strict lower bounds are placed on 
sectional curvature relative to the mean curvature. For 
further details, see \cite{ER, GUXU, xuzhao}.

The aim of this paper is to explore the topology and extrinsic 
geometry of submanifolds under a lower bound on sectional 
curvature expressed in terms of the mean curvature, without 
imposing any conditions on the mean curvature. In contrast 
to all previously known results, the pinching condition on the 
sectional curvature is no longer required to be strict.

For any Riemannian manifold $M^n$, we define the function 
$K_{\text{min}}$ as 
$$
K_{\text{min}}(x)=
\min\left\{K(\pi)\colon \pi \,\,{\text {is a 2-plane in}}\,\,T_xM\right\},
$$
where $K(\pi)$ denotes the sectional curvature of the plane 
$\pi$ at the point $x\in M^n$.
Consider an isometric immersion $f\colon M^n\to\Q_c^{n+m}$,
where $\Q_c^{n+m}$ is a space form of constant curvature 
$c$. 
The second fundamental form of $f$, which takes values 
in the normal bundle, is denoted by 
$\a_f\colon TM\times TM\to N_fM$. The mean 
curvature vector field of $f$ is defined by 
$\mathcal H=(1/n)\,\trace\a_f$, where 
$\mathrm{tr}$ denotes the trace. The corresponding 
mean curvature function is defined by 
$H=\|\mathcal H\|$. Throughout the paper, all 
submanifolds under consideration are assumed 
to be connected. Our result can now be stated as follows.

\begin{theorem}\label{thm1} 
Let $f\colon M^n\to\Q_c^{n+m}$, 
$n\geq 4$, $c\geq 0$, be an isometric immersion of a compact, 
oriented Riemannian manifold. Assume that, at every point of 
$M^n$, the sectional curvature satisfies 
\be\label{*}
K_{\text{min}}\geq b(n,H,c)=
\begin{cases}
\,\frac{1}{2}(c+\frac{n}{4}H^2)&\text{if\, } n\geq 5 \\[1mm]
\,\frac{1}{3}(c+H^2)&\text{if\, } n=4.
\end{cases}
\tag{$\ast$}
\ee 
Then the following assertions hold:
\vspace{1ex} 

\noindent $(i)$ If $n\geq 5$, then $M^n$ is homeomorphic to 
the sphere $\Sf^n$. 
\vspace{1ex}

\noindent $(ii)$ If $n=4$, then either $M^4$ is diffeomorphic 
to $\Sf^4$, or equality holds in \eqref{*} at every point, and 
$M^4$ is isometric to the projective plane $\CP_r^2$ 
of constant holomorphic curvature $\frac{4}{3r^2}$ with 
$r=\frac{1}{\sqrt{c+H^2}}$. Furthermore, the immersion $f$ 
is a composition $f=j\circ g$, where $g$ is 
the standard embedding of $\CP_r^2$ into a 
sphere $\Sf^7(r)$ and $j\colon\Sf^7(r)\to\Q_c^{4+m}$ 
is an umbilical inclusion. 
\end{theorem}

It is noteworthy that our pinching condition \eqref{*} 
fully determines the geometry in the case of 
four-dimensional submanifolds that are not 
diffeomorphic to $\Sf^n$, as described in part 
$(ii)$ of the above theorem. In contrast, 
for all other dimensions, the condition \eqref{*} 
allows us to describe only the topology of the 
submanifold. Interestingly, striking differences 
emerge not only in the results but also in the 
arguments between the two cases $n = 4$ and 
$n \geq 5$.

For dimensions $n\geq5$, the strategy relies 
on seminal results of Lawson and Simons \cite{LS} 
which addressed compact submanifolds in space forms 
of nonnegative curvature. They 
proved that under certain upper 
bounds of the second fundamental 
form of the submanifold, specific integral homology 
groups must vanish. We show that our pinching 
condition ensures the validity of the bounds 
established by Lawson and Simons. 
In particular, for dimensions $n=5,6,12,56,61$, the 
manifold $M^n$ is not only topologically equivalent to 
a sphere but also diffeomorphic to it. This follows from 
the fact that, in these cases, the differentiable structure 
of the sphere is known to be unique, as stated in 
Corollary $1.15$ of \cite{WX}.

In the case of four-dimensional submanifolds, the proof of 
our result is more elaborate. It heavily relies on concepts 
from four-dimensional geometry, the geometry of Riemannian 
manifolds with nonnegative isotropic curvature, and the 
Bochner technique, all of which play a crucial role.

The Veronese submanifold (cf. \cite{Itoh1, Itoh2}) is a minimal 
isometric embedding $\psi_n\colon\R P^n_K\to\Sf^{n+m}$ 
of the real projective space $\R P^n_K$ of sectional 
curvature $K=\frac{n}{2(n+1)}$ into a unit sphere. Clearly, this 
immersion satisfies \eqref{*} if $n=4$, illustrating that 
the assumption of the manifold being oriented is essential 
and cannot be omitted. In the final section, we present a 
method for constructing geometrically distinct 
submanifolds that are topological spheres and 
satisfy the condition \eqref{*}. These examples 
show that, even under the strict pinching condition, 
the spherical conclusion cannot in general be strengthened 
to an extrinsic rigidity statement.

It is also worth mentioning that the geometric and 
topological rigidity of submanifolds has been studied 
recently \cite{DV2, ge} under lower bounds on the Ricci 
curvature in terms of the mean curvature, without requiring 
strict inequalities. Notably, there exist submanifolds that 
satisfy our pinching condition on the sectional curvature 
but fail to meet the pinching conditions on the Ricci curvature 
in \cite{DV2, ge}. This demonstrates that our result does 
not follow from those in \cite{DV2, ge}.

\section{The case $n\geq5$}

The result by Lawson and Simons \cite{LS} aforementioned
was strengthened by Elworthy and Rosenberg \cite[p.\ 71]{ER} 
by not requiring the bound to be strict at all points of 
the submanifold. In this section, we first state their result
and then analyse the relation between the condition required 
in their theorems and our pinching assumption.

\begin{theorem}{\em(\cite{ER, LS, X})}\label{ls} 
Let $f\colon M^n\to\Q_c^{n+m}$, $c\geq 0$, 
be an isometric immersion of a compact manifold and 
let $p$ be an integer such that $1\leq p\leq n-1$. 
Assume that at every point $x\in M^n$ and for 
every orthonormal basis $\{e_1,\dots,e_n\}$ of 
$T_xM$, the second fundamental form 
$\alpha_f\colon TM\times TM\to N_fM$ satisfies
\be\label{LS}
\Theta_p=\sum_{i=1}^p\sum_{j=p+1}^n\big(2\|\a_f(e_i,e_j)\|^2
-\<\a_f(e_i,e_i),\a_f(e_j,e_j)\>\big)\leq p(n-p)c \tag{$**$}.
\ee
If at a point of $M^n$ the inequality \eqref{LS} is strict for
every orthonormal basis, then the homology groups satisfy 
$H_p(M^n;\mathbb{Z})=H_{n-p}(M^n;\mathbb{Z})=0$.
\end{theorem}

The following result is a key ingredient in the 
proof of our theorem in higher dimensions.

\begin{proposition}\label{prop}
Let $f\colon M^n\to\Q^{n+m}_c$, $n\geq 5$, $c\geq 0$, 
be an isometric immersion satisfying the inequality 
\eqref{*} at a point $x\in M^n$. Then, the following 
assertions hold at $x$:
\vspace{1ex}

\noindent $(i)$ The inequality \eqref{LS} is satisfied for all 
$2\leq p\leq n-2$ and every orthonormal 
basis $\{e_i\}_{1\leq i\leq n}$ of $T_xM$. 
Moreover, if \eqref{*} is strict at $x$, then 
the inequality \eqref{LS} is strict for all $2\leq p\leq n-2$ and 
every orthonormal basis of $T_xM$. 
\vspace{1ex}

\noindent $(ii)$ If equality holds in \eqref{LS} for some 
orthonormal basis $\{e_i\}_{1\leq i\leq n}$ of 
$T_xM$ and for an integer $2\leq p\leq n-2$, 
then $H(x)=0$, and the sectional curvature 
satisfies
$$
K(e_i\wedge e_j)=K_{\text{min}}(x)=\frac{c}{2},\quad {\text {for all}}\,\,\,\,1\leq i\leq p \,\,{\text {and}}\,\,p+1\leq j\leq n. 
$$
\end{proposition}

\proof
Let $\{e_i\}_{1\leq i\leq n}$ be an arbitrary orthonormal basis 
of $T_xM$. For convenience, write $\a_{ij}=\alpha_f(e_i,e_j)$, 
$1\leq i,j\leq n$. By using the Gauss equation (cf. \cite{DT}), 
we then have
\begin{align}\label{a}
\Theta_p&=\;\sum_{i=1}^p\sum_{j=p+1}^n\big(2c-2K(e_i\wedge e_j)
+\<\a_{ii},\a_{jj}\>\big)\nonumber\\
&\leq\; 2p(n-p)\big(c-K_{\text{min}}(x)\big)-
\|\sum_{i=1}^p\a_{ii}-\frac{n}{2}\mathcal H_f(x)\|^2
+\frac{n^2}{4}H^2(x).
\end{align}
Since $p(n-p)>n$ for all $2\leq p\leq n-2$, it follows from \eqref{a} and our 
assumption on the sectional curvature that 
\begin{align}\label{a1}
\Theta_p-p(n-p)c 
\leq\;p(n-p)\Big(c+\frac{n}{4}H^2(x)-2K_{\text{min}}(x)\Big)\leq0.
\end{align}

If inequality \eqref{*} is strict, then also \eqref{a} and 
thus \eqref{a1} becomes 
strict, and this completes the proof of part $(i)$.

Next, we assume that equality holds in \eqref{LS} for 
some orthonormal basis $\{e_i\}_{1\leq i\leq n}$ of 
$T_xM$ and an integer $2\leq p\leq n-2$. In this case, 
equality holds in both \eqref{a} and \eqref{a1}. 
Since $p(n-p)>n$, it follows that $H(x)=0$, and 
$$
K(e_i\wedge e_j)=K_{\text{min}}(x)=\frac{c}{2},\quad{\text {for all}}\,\,\,
1\leq i\leq p \,\,{\text {and}}\,\,p+1\leq j\leq n,
$$
which completes the proof of part $(ii)$.\qed
\vspace{1ex}

\begin{lemma}\label{n>4}
Let $f\colon M^n\to\Q^{n+m}_c$, $n\geq 4$, $c\geq 0$, 
be an isometric immersion of a compact Riemannian 
manifold satisfying the inequality \eqref{*}. Then the 
following assertions hold:
\vspace{1ex}

\noindent $(i)$The universal cover of $M^n$ is compact.
\vspace{1ex}

\noindent $(ii)$ If $M^n$ is oriented and $H_p(M^n,\Z)=0$ 
for all $2\leq p\leq n-2$, then $M^n$ is a homology sphere. 
\end{lemma}
\proof \noindent $(i)$ First, assume that $c>0$. Then, our 
assumption implies that the sectional curvature of $M^n$ 
is bounded below by a positive constant. By the 
Bonnet-Myers theorem, the universal cover of $M^n$ 
is compact. 

Now, assume that $c=0$. Clearly, under our assumptions, 
the manifold $M^n$ has nonnegative Ricci curvature. 
According to the Cheeger-Gromoll splitting theorem 
\cite{cg}, the universal cover of $M^n$, equipped with 
the covering metric, splits isometrically as a product 
$\R^p\times N$, where $N$ is a compact manifold. 
Furthermore, there exists a point $x\in M^n$ such 
that $H(x)>0$, since there are no compact minimal 
submanifolds in Euclidean space. This implies that 
the Ricci curvature is positive in every tangent 
direction at $x$. Hence, the universal cover cannot 
contain a nontrivial Euclidean factor. Therefore, the 
universal cover of $M^n$ is compact.

\indent $(ii)$ Suppose now that $M^n$ is oriented and 
$H_p(M^n,\Z)=0$ for all $2\leq p\leq n-2$. By the 
universal coefficient theorem of cohomology, 
$H^{n-1}(M^n;\Z)$ is torsion-free. Using Poincar\'e 
duality, it follows that $H_1(M^n;\Z)$ is also 
torsion-free. The compactness of the universal cover 
of $M^n$ implies that the group of deck transformations, 
and hence the fundamental group $\pi _1(M^n)$ is
finite. Consequently, $H_1(M^n;\Z)=0$. Since $M^n$ 
is oriented, we have $H^n(M^n;\Z)\cong\Z$. 
Applying the universal coefficient theorem, we deduce 
that $H_{n-1}(M^n;\Z)$ is torsion-free. By Poincar\'e 
duality, we conclude that $H_{n-1}(M^n;\Z)=0$, 
which implies that $M^n$ is a homology 
sphere.\qed
\vspace{1ex}

\noindent\emph{Proof of Theorem \ref{thm1} for $n\geq5$.} 
We first claim that $H_p(M^n,\Z)=0$ for all $2\leq p\leq n-2$. 
Suppose to the contrary that $H_p(M^n,\Z)\neq0$ for 
some integer $2\leq p\leq n-2$. From Proposition 
\ref{prop}, we know that the inequality \eqref{LS} is 
satisfied for any integer $2\leq p\leq n-2$ and for 
every orthonormal basis of 
$T_xM$ at every point $x\in M^n$. Since 
$H_p(M^n,\Z)\neq0$, it follows from 
Theorem \ref{ls} that at every point there exists an 
orthonormal basis such that inequality \eqref{LS} 
holds as equality. By part $(ii)$ of Proposition 
\ref{prop}, this implies that the submanifold $f$ 
is minimal, and hence $c>0$, with sectional curvature 
$$
K\geq \frac{c}{2}>\frac{nc}{2(n+1)}.
$$ 
From a well known result 
(see, e.g., \cite{Itoh1, Itoh2} or \cite{Leung}), it follows that 
the submanifold is totally geodesic. Consequently, 
$M^n$ is isometric to the round sphere, which is a 
contradiction. Thus, we conclude that $H_p(M^n,\Z)=0$ 
for all $2\leq p\leq n-2$. Hence, by Lemma \ref{n>4}, 
$M^n$ is a homology sphere.

After passing to the universal covering 
$\pi\colon\tilde M^n\to M^n$, and noting that the 
isometric immersion $\tilde f=f\circ\pi$ satisfies 
the inequality \eqref{*}, we proceed as follows. 
Since by Lemma \ref{n>4} the universal cover 
$\tilde M^n$ is compact, we argue as above 
and conclude that $\tilde M^n$ is a homology sphere. 
Thus $\tilde M^n$ is a simply connected homotopy 
sphere. By the resolution of Poincar\'e conjecture, 
we know that $\tilde M^n$ is homeomorphic to 
$\Sf^n$. Since the homology sphere $M^n$ is 
covered by $\Sf^n$, a result by Sjerve 
\cite{Sj} implies that $\pi _1(M^n)=0$, and 
consequently, $M^n$ is homeomorphic 
to $\Sf^n$.\qed

\section{Background}

\subsection{Geometry of four-dimensional manifolds}
In this section, we collect basic facts about four-dimensional 
geometry. For a detailed exposition of the subject, we refer 
the reader to \cite{LeB, sea}.

Let $(M,\<\cdot,\cdot\>)$ be an oriented Riemannian 
manifold of dimension $n=4$ with Levi-Civita connection 
$\n$ and curvature tensor $R$ given by
$$
R(X,Y)=[\n_X,\n_Y]-\n_{[X,Y]},\quad X,Y\in\mathcal X(M).
$$
The Ricci tensor of $(M,\<\cdot,\cdot\>)$ is defined by
$$
\Ric(X,Y)=\sum_i\<R(X,E_i)E_i,Y\>,\quad X,Y\in\mathcal X(M),
$$
where $\{E_i\}_{1\leq i\leq 4}$ is a local orthonormal frame.

At every point $x\in M$, we consider the \emph {Bochner-Weitzenb\"ock 
operator} $\B^{[2]}$ as an endomorphism of the space of 
2-vectors $\Lambda^2T_xM$ at $x$ given by 
\begin{align}\label{WB}
\<\<\B^{[2]}(v_1\wedge v_2), w_1\wedge w_2\>\>
=&\;\Ric(v_1,w_1)\<v_2,w_2\>+\Ric(v_2,w_2)\<v_1,w_1\>\nonumber\\
&-\Ric(v_1,w_2)\<v_2,w_1\>-\Ric(v_2,w_1)\<v_1,w_2\>\nonumber\\
&-2\<R(v_1,v_2)w_2,w_1\>,
\end{align}
and then extend it linearly to all of $\Lambda^2T_xM$.
Here $\<\<\cdot,\cdot\>\>$ stands for the inner product of 
$\Lambda^2T_xM$ defined by 
$$
\<\<v_1\wedge v_2,w_1\wedge w_2\>\>=\det(\<v_i,w_j\>).
$$
The Bochner-Weitzenb\"ock operator is a self-adjoint operator. 
If $X\in\Lambda^2T_xM$, the dual 2-form $\omega$ is 
defined by $\omega(v,w)=\<\<X,v\wedge w\>\>$, 
and we may consider $X$ as the skew-symmetric 
endomorphism of the tangent space at $x$ by 
$$
\<X(v),w\>=\omega(v,w)=\<\<X,v\wedge w\>\>. 
$$

Clearly $\B^{[2]}$ can also be viewed as an endomorphism 
of the bundle $\Omega^2(M)$ of 2-forms of the manifold via 
$\<\<\cdot,\cdot\>\>$. If $\omega$ is a 2-form, then 
$\B^{[2]}(\omega)$ is given by 
$$
\B^{[2]}(\omega)(X_1,X_2)=\omega(\Ric(X_1),X_2)
+\omega(X_1,\Ric(X_2))-\sum_i\omega(R(X_1,X_2)E_i,E_i).
$$
Then the Bochner-Weitzenb\"ock operator acts on 
$\Lambda^2TM$ by 
$$
\omega(\B^{[2]}(X_1\wedge X_2))
=\B^{[2]}(\omega)(X_1,X_2).
$$ 
Taking $\omega$ to 
be the dual form to $w_1\wedge w_2$ yields \eqref{WB}. 
Throughout the paper, we will mostly identify 2-forms 
with 2-vectors. 

We recall the Bochner-Weitzenb\"ock formula which 
can be written as 
$$
\<\Delta\omega,\omega\>=\frac{1}{2}\Delta\|\omega\|^2+\|\n\omega\|^2
+\<\B^{[2]}(\omega),\omega\>
$$
for every $\omega\in\Omega^2(M)$. This implies that any 
harmonic 2-form on a compact manifold is parallel provided 
that the Bochner-Weitzenb\"ock operator is nonnegative. 

The bundle of 2-forms of any oriented four-dimensional 
Riemannian manifold $M$ decomposes as a direct sum 
$$
\Omega^2(M)=\Omega_+^2(M)\oplus\Omega_-^2(M)
$$
of the eigenspaces of the Hodge star operator 
$\ast\colon\Omega^2(M)\to\Omega^2(M)$. The sections of 
$\Omega_+^2(M)$ are called {\emph {self-dual 2-forms}}, 
whereas the ones of $\Omega_-^2(M)$ are called 
{\emph {anti-self-dual 2-forms}}. 
Accordingly, we have the splitting
$$
\Lambda^2T_xM=\Lambda_+^2T_xM\oplus\Lambda_-^2T_xM
$$
at every point $x$, where $\Lambda_\pm^2T_xM$ 
are the eigenspaces of the Hodge star operator 
$\ast\colon\Lambda^2T_xM\to\Lambda^2T_xM$. Both 
spaces $\Lambda_+^2T_xM$ and $\Lambda_-^2T_xM$ are 
$\B^{[2]}$-invariant (see Proposition 1 in \cite{smichig}). Then we have 
accordingly the decomposition $\B^{[2]}=\B_+^{[2]}\oplus\B_-^{[2]}$.

Suppose now that $M$ is a compact 
oriented Riemannian four-manifold. The Hodge theorem 
guarantees that every de Rham cohomology class on $M$ has a unique 
harmonic representative. In particular, the space 
$\mathscr H^2(M)$ of harmonic 2-forms decomposes as 
$$
\mathscr H^2(M)=\mathscr H^2_+(M)\oplus\mathscr H^2_-(M),
$$ 
where $\mathscr H^2_+(M)$ and $\mathscr H^2_-(M)$ are the spaces of 
self-dual and anti-self-dual harmonic 2-forms, respectively. The dimensions 
of these subspaces, denoted by
$\beta_\pm(M)=\dim\mathscr H^2_\pm(M)$, are oriented 
homotopy invariants of $M$. Their difference 
$\tau=\beta_+(M)-\beta_-(M)$ is the signature of $M$, 
while their sum equals the second Betti number 
$\beta_2(M)$ of the manifold $M$. 

\subsection{Isotropic curvature}
Let $(M,g)$ be a Riemannian manifold of dimension $n\geq4$. 
We say that $(M,g)$ has \emph {nonnegative isotropic curvature} 
at a point $x\in M$ if
$$
R_{1331} + R_{1441}+ R_{2332} + R_{2442}
-2 R_{1234} \geq0, 
$$
for all orthonormal four-frames 
$\{e_1,e_2,e_3,e_4\}\subset T_xM$. Here, we denote 
$$
R_{ijk\ell}=g\big(R(e_i,e_j)e_k, e_\ell\big), \quad1\leq i,j,k,\ell\leq4, 
$$ where $R$ is the 
curvature tensor. If the 
strict inequality holds, we say that $(M,g)$ has 
\emph{positive isotropic curvature} at $x$. 
The manifold $(M,g)$ is said to have nonnegative (or positive) 
isotropic curvature if it satisfies the corresponding condition 
at every point and for all orthonormal four-frames.

The following result is well known (see for instance \cite[p. 161]{Nor} 
or Section 2 of \cite{MWolf}), and its proof is implicitly contained in the 
argument of Proposition 2.1 in \cite{MWolf}.

\begin{lemma}\label{sb}
Let $M$ be a four-dimensional Riemannian manifold. Then, 
at every point $x \in M$, the nonnegativity of the isotropic 
curvature is equivalent to the nonnegativity of the 
Bochner-Weitzenb\"ock operator $\B^{[2]}$.
\end{lemma}

Proposition 1.1 of \cite{sea} also establishes, in a more 
general setting, that nonnegative isotropic curvature 
implies the nonnegativity of the Bochner-Weitzenb\"ock 
operator $\B^{[2]}$. This is precisely the direction 
needed for our proofs.

\section{Proof of Theorem \ref{thm1} for the case $n=4$}

We need the following auxiliary lemma.
\begin{lemma}\label{R}
Let $M^4$ be a four-dimensional Riemannian manifold. 
Suppose that there exists an orthonormal four-frame 
$\{e_1,e_2,e_3,e_4\}\subset T_pM$ at a point 
$p\in M^4$ such that the sectional curvature
$K_{ij}=K(e_i\wedge e_j)$ satisfies 
$$
K_{ij}=K_{\text{min}}(p),\quad
\text{for}\;\; i=1,2 \;\;\text{and}\;\; j=3,4.
$$ 
Then the curvature tensor satisfies $R_{ijk\ell}=0$ 
whenever three of the indices $i,j,k,\ell$ 
are mutually distinct. 
\end{lemma}
\proof 
We consider the function $F\colon U\to\R$ defined by 
\begin{align*}
F(x,y)&=\<R\big(\tau(x),\tau(y)\big)\tau(y),\tau(x)\>-K_{\text{min}}(p)\|\tau(x)\wedge\tau(y)\|^2\\
&=\sum_{i,j,k,\ell=1}^4x_ix_jy_ky_\ell R_{ik\ell j}-K_{\text{min}}(p)\big(\|x\|^2\|y\|^2-\<x,y\>^2\big), 
\end{align*}
where 
$$
U=\left\{(x,y)\in\R^8: x\wedge y\neq0\right\},\quad
x=(x_1,x_2,x_3,x_4),\quad y=(y_1,y_2,y_3,y_4),
$$ 
and $\tau\colon\R^4\to T_pM$ 
is the linear isometry given by $\tau(x)=\sum_{i}x_ie_i$.

By our assumption, the nonnegative function $F$ attains its minimum 
at $\varepsilon_{ij}=(\tau^{-1}e_i,\tau^{-1}e_j)\in U$, for $i=1,2$ and $j=3,4$. 
This implies that
\be\label{fermat}
\frac{\partial F}{\partial x_k}(\varepsilon_{ij})
=\frac{\partial F}{\partial y_k}(\varepsilon_{ij})=0,\quad\text{for}\;\;i=1,2,\;
j=3,4\;\;\text{and}\;\;1\leq k\leq4. 
\ee
Using the algebraic properties of the curvature tensor, a direct 
computation yields 
$$
\frac{\partial F}{\partial x_k}(\varepsilon_{ij})
=2R_{kjji}\;\,\text{and}\;\,\frac{\partial F}{\partial y_k}(\varepsilon_{ij})
=2R_{kiij},\quad\text{for}\;\;i=1,2,\;j=3,4\;\;\text{and}\;\, k\neq i,j.
$$
From this, along with \eqref{fermat}, and the algebraic 
properties of the curvature tensor, it follows that 
$R_{ijk\ell}=0$ whenever three of the indices 
$i,j,k,\ell$ are mutually distinct.\qed
\vspace{1ex}

The following result is essential for the proof of the theorem for $n=4$. 

\begin{proposition}\label{nnic}
Let $f\colon M^4\to\Q^{4+m}_c$, $c\geq 0$, be an 
isometric immersion of an oriented four-dimensional 
Riemannian manifold satisfying the inequality 
\eqref{*} at a point $x\in M^4$. Then the 
following assertions hold:
\vspace{1ex}

\noindent $(i)$ The manifold $M^4$ has nonnegative 
isotropic curvature at $x$.
\vspace{1ex}

\noindent $(ii)$ Suppose that $M^4$ has 
not positive isotropic curvature at $x$. Then 
equality holds in \eqref{*} at $x$, 
and there exists an oriented orthonormal four-frame 
$\{e_1,e_2,e_3,e_4\}\subset T_xM$ such that 
the following assertions hold at $x$:
\begin{itemize}
\item[(ii1)] The sectional curvature satisfies 
\be\label{K12}
K_{ij}=K_{\text{min}}(x),\quad
\text{for}\;\; i=1,2 \;\;\text{and}\;\; j=3,4.
\ee
\item[(ii2)] The second fundamental form 
$\a_f$ of $f$ satisfies the following conditions: 
\begin{align}
\a_{11}=\a_{22},\;\,\a_{33}=\a_{44},\;\,\a_{12}&=\a_{34}=0,\;\,\a_{23}
=-\varepsilon\a_{14},\;\,\a_{24}=\varepsilon\a_{13},\label{B2a1e}\\
\|\a_{13}\|=\|\a_{14}\|,\;\,\<\a_{13},\;\a_{14}\>&=\<\a_{jj},\a_{1i}\>=0,\quad 
1\leq j\leq4,\;i=3,4,\label{B2a2e}\\
4\big(c+\<\a_{11},a_{44}\>\big)&=\frac{1}{2}\,\|\a_{11}-a_{44}\|^2+6\|\a_{13}\|^2,\label{B2a3}
\end{align}
where, for convenience, we write $\a_{ij}=\a_f(e_i,e_j)$ and 
$\varepsilon=\pm1$.
\item[(ii3)] Consider the orthonormal basis 
$\{\eta_i\}_{1\leq i\leq6}$ of the space 
of 2-vectors $\Lambda^2T_xM$, satisfying 
$\eta_i\in\Lambda_+^2T_xM$ and 
$\eta_{i+3}\in\Lambda_-^2T_xM$ for $1\leq i\leq 3$, defined by 
\begin{align*}
\eta_1=\frac{1}{\sqrt{2}}(e_{12}+e_{34}),\quad
\eta_2&=\frac{1}{\sqrt{2}}(e_{13}-e_{24}),
\quad\eta_3=\frac{1}{\sqrt{2}}(e_{14}+e_{23}),\\
\eta_4=\frac{1}{\sqrt{2}}(e_{12}-e_{34}),\quad
\eta_5&=\frac{1}{\sqrt{2}}(e_{13}+e_{24}),\quad
\eta_6=\frac{1}{\sqrt{2}}(e_{14}-e_{23}),
\end{align*}
where $e_{ij}=e_i\wedge e_j$. This basis diagonalizes the 
Bochner-Weitzenb\"ock operator $\B^{[2]}$ 
at the point $x$, with corresponding 
eigenvalues $\{\mu_i\}_{1\leq i\leq6}$ given by: 
\begin{align*}
\mu_1=\;&4\big(K_{\text{min}}(x)-\varepsilon\|\a_{13}\|^2\big)\geq0,\\
\mu_2=\;&\mu_3=K_{12}+K_{34}+2\big(K_{\text{min}}(x)+\varepsilon\|\a_{13}\|^2\big)\geq0,\\
\mu_4=\;&4\big(K_{\text{min}}(x)+\varepsilon\|\a_{13}\|^2\big)\geq0,\\
\mu_5=\;&\mu_6=K_{12}+K_{34}+2\big(K_{\text{min}}(x)-\varepsilon\|\a_{13}\|^2\big)\geq0.
\end{align*}
\end{itemize}
\end{proposition}
\proof 
\noindent $(i)$ Clearly, the manifold $M^4$ has nonnegative 
isotropic curvature at $x$ if 
$$
R_{1331} + R_{1441}+ R_{2332} + R_{2442}
+2 R_{1234} \geq0 
$$
for every orthonormal four-frame at $x$. 
Let $\{e_1,e_2,e_3,e_4\}\subset T_xM$ 
be an arbitrary orthonormal four-frame, and let 
$\{\xi_\a\}_{1\leq \a\leq m}$ be an orthonormal basis 
of the normal space $N_fM(x)$ of $f$ at $x$. 
For simplicity, we set $h^\a_{ij}=\<\a_{ij},\xi_\a\>$. 
Using the Gauss equation, we have 
\begin{align}\label{ic1}
R_{1331} + R_{2332}-&R_{1243}
\geq 2K_{\text{min}}(x)+\sum_\a(h^\a_{14}h^\a_{23}-h^\a_{13}h^\a_{24})
\nonumber\\
&\geq 2K_{\text{min}}(x)
-\frac{1}{2}\sum_\a\big((h^\a_{13})^2+(h^\a_{23})^2+(h^\a_{14})^2+(h^\a_{24})^2\big). 
\end{align}

Following the computations in \cite[p.~341]{GUXU}, we obtain
\begin{align*}
16K_{\min}(x)&\le 4K\big((e_1+e_3)\wedge(e_2-e_4)\big)
+4K\big((e_1-e_3)\wedge(e_2+e_4)\big)\\
&\quad +4K\big((e_1+e_4)\wedge(e_2+e_3)\big)
+4K\big((e_1-e_4)\wedge(e_2-e_3)\big)\\
&=2(R_{1331}+R_{1441}+R_{2332}+R_{2442})
+4(R_{1221}+R_{3443})-12R_{1243}.
\end{align*}
Consequently,
\be\label{ic2}
6R_{1243}\leq R_{1331}+R_{1441}+R_{2332}+R_{2442}
+2(R_{1221}+R_{3443})-8K_{\text{min}}(x).\nonumber
\ee
From this inequality, we deduce 
\begin{align}\label{ic3t}
R_{1331}&+ R_{2332}-R_{1243}
\geq\frac{10}{3}K_{\text{min}}(x)\nonumber\\
&-\frac{1}{6}(R_{1331}+R_{1441}+R_{2332}+R_{2442}+2R_{1221}+2R_{3443}).
\end{align}
Using the Gauss equation, it follows from \eqref{ic3t} that 
\begin{align}\label{ic3}
R_{1331} &+ R_{2332}-R_{1243}
\geq\frac{10}{3}K_{\text{min}}(x)-\frac{4}{3}c+\frac{1}{3}\sum_\a\big((h^\a_{12})^2+(h^\a_{34})^2\big)\nonumber\\
&-\frac{1}{6}\sum_\a\big(h^\a_{11}h^\a_{33}-(h^\a_{13})^2+h^\a_{11}h^\a_{44}-(h^\a_{14})^2
\big)\nonumber\\
&-\frac{1}{6}\sum_\a\Big(h^\a_{22}h^\a_{33}-(h^\a_{23})^2+h^\a_{22}h^\a_{44}-(h^\a_{24})^2
+2h^\a_{11}h^\a_{22}+2h^\a_{33}h^\a_{44}\Big)\nonumber\\
&\geq \frac{10}{3}K_{\text{min}}(x)-\frac{4}{3}c
+\frac{1}{6}\sum_\a\big((h^\a_{13})^2+(h^\a_{23})^2+(h^\a_{14})^2+(h^\a_{24})^2\big)
\nonumber\\
&-\frac{1}{6}\sum_\a\big((h^\a_{11}+h^\a_{33})(h^\a_{22}+h^\a_{44})+(h^\a_{11}+h^\a_{44})(h^\a_{22}+h^\a_{33})
\big). 
\end{align}
From $\frac{1}{4}\times$\eqref{ic1}+$\frac{3}{4}\times$\eqref{ic3}, we obtain 
\begin{align}\label{ic4}
R_{1331} + R_{2332}&-R_{1243}
\geq3K_{\text{min}}(x)-c\nonumber\\
&-\frac{1}{8}\sum_\a\big((h^\a_{11}+h^\a_{33})(h^\a_{22}+h^\a_{44})+(h^\a_{11}+h^\a_{44})(h^\a_{22}+h^\a_{33})
\big). 
\end{align}

On the other hand, we have 
\be\label{ic5}
(h^\a_{11}+h^\a_{33})(h^\a_{22}+h^\a_{44})
=-\big(h^\a_{11}+h^\a_{33}-\frac{1}{2}\trace A_\a\big)^2
+\frac{1}{4}(\trace A_\a)^2
\leq\frac{1}{4}(\trace A_\a)^2
\ee
and similarly
\be\label{ic6}
(h^\a_{11}+h^\a_{44})(h^\a_{22}+h^\a_{33})
=-\big(h^\a_{11}+h^\a_{44}-\frac{1}{2}\trace A_\a\big)^2
+\frac{1}{4}(\trace A_\a)^2\leq\frac{1}{4}(\trace A_\a)^2, 
\ee
where $A_\a$ denotes the shape operator associated 
to $\xi_\a$. If both \eqref{ic5} and \eqref{ic6} hold as 
equalities, then 
$$
h^\a_{11}+h^\a_{33}=\frac{1}{2}\trace A_\a=h^\a_{11}+h^\a_{44}, 
$$
from which we deduce that $h^\a_{11}=h^\a_{22}$ and 
$h^\a_{33}=h^\a_{44}$. 

Then using \eqref{ic4}, \eqref{ic5}, and \eqref{ic6}, we obtain 
\be\label{ic7}
R_{1331} + R_{2332}-R_{1243}
\geq3K_{\text{min}}(x)-(c+H^2(x)).
\ee
Observe that if equality holds in \eqref{ic7}, then all 
inequalities \eqref{ic1}--\eqref{ic6} also become equalities. 
Hence, we have 
\be\label{ic3q11}
R_{1331}=R_{2332}=K_{\text{min}}(x)
\ee
and
\begin{align}\label{ic3q1}
 h^\a_{14}+h^\a_{23}=0,\,\,\,\,h^\a_{13}&=h^\a_{24},\,\,\,\,h^\a_{12}=h^\a_{34}=0,
\,\,\,\,h^\a_{11}=h^\a_{22},\,\,\,\,h^\a_{33}=h^\a_{44},
\end{align}
for all $1\leq\a\leq m$. 

In a similar manner, by working instead with the four-frame 
$\{e_1,e_2,e_4,-e_3\}$ we obtain 
\be\label{ic8}
R_{1441} + R_{2442}-R_{1243}
\geq3K_{\text{min}}(x)-(c+H^2(x))
\ee
and if equality holds in the above inequality, then \eqref{ic3q1} 
holds and 
\be\label{ic3q12}
R_{1441}=R_{2442}=K_{\text{min}}(x).
\ee

From \eqref{ic7}, \eqref{ic8} and \eqref{*}, it follows that 
\be\label{ic}
R_{1331} + R_{1441}+R_{2332} + R_{2442}-2R_{1243}
\geq2\big(3K_{\text{min}}(x)-(c+H^2(x))\big)\geq0.
\ee
This completes the proof of part $(i)$. 
Note that if inequality \eqref{*} is strict, then \eqref{ic} 
implies that the manifold $M^4$ has positive 
isotropic curvature at $x$.
\vspace{1ex}

\indent $(ii)$ Suppose now that the manifold 
$M^4$ does not have positive isotropic 
curvature at $x$. Then equality holds in 
\eqref{*} at $x$, and there exists an 
orthonormal four-frame 
$\{e_1,e_2,e_3,e_4\}\subset T_xM$ such 
that equality holds in \eqref{ic}. Hence, both 
\eqref{ic7} and \eqref{ic8} are equalities, 
which implies that \eqref{ic3q11}, 
\eqref{ic3q1}, and \eqref{ic3q12} hold. 
Consequently, \eqref{K12} is satisfied, 
completing the proof of part $(ii1)$. 

Moreover, \eqref{ic3q1} implies 
that the second fundamental 
form satisfies 
\be\label{B2a1}
\a_{11}=\a_{22},\quad\a_{33}=\a_{44},\quad\a_{12}=
\a_{34}=0,\quad\a_{13}=\a_{24},\quad\a_{14}+\a_{23}=0.
\ee
From \eqref{K12}, we have $K_{13}=K_{14}$. Together 
with the second equality in \eqref{B2a1}, the Gauss 
equation then yields 
\be\label{B2a11}
\|\a_{13}\|=\|\a_{14}\|.
\ee 
Next, \eqref{K12} and Lemma \ref{R} imply that $R_{ijk\ell}=0$ whenever three of the indices $i,j,k,\ell$ are mutually distinct. Since $R_{1223}=R_{1443}=0$, \eqref{B2a1} and the Gauss equation yield 
\be\label{B2a12}
\<\a_{13},\a_{11}\>=\<\a_{13},\a_{44}\>=0.
\ee 
Similarly, from $R_{1224}=R_{1334}=0$,
we obtain 
\be\label{B2a13}\<\a_{14},\a_{11}\>=\<\a_{14},\a_{44}\>=0.
\ee 
Finally, $R_{3114}=0$ implies that 
\be\label{B2a14}\<\a_{13},\a_{14}\>=0.
\ee 

From \eqref{K12}, the Gauss equation, and the second equality 
in \eqref{B2a1}, we obtain 
$$
K_{\text{min}}(x)=K_{13}=c+\<\a_{11},\a_{44}\>-\|\a_{13}\|^2.
$$
Moreover, the first two equalities in 
\eqref{B2a1} yield 
$$
4H^2=\|\a_{11}-\a_{44}\|^2+4\<\a_{11},\a_{44}\>.
$$
Equation \eqref{B2a3} follows immediately from the 
preceding two equalities, since \eqref{*} holds as 
an equality at $x$.

If the four-frame $\{e_1,e_2,e_3,e_4\}$ 
is positively oriented, 
then \eqref{B2a1e} and \eqref{B2a2e} 
follow immediately from \eqref{B2a1}--\eqref{B2a14} 
with $\varepsilon=1$. Otherwise, we 
replace the frame by $\{e_1,e_2,e_3,-e_4\}$. 
Clearly, this frame satisfies \eqref{B2a1e}--\eqref{B2a3} 
with $\varepsilon=-1$, thereby completing the proof of $(ii2)$. 

By a straightforward computation using \eqref{WB}, 
the Gauss equation, \eqref{B2a1e}, and \eqref{B2a2e}, 
we obtain that the Bochner-Weitzenb\"ock operator 
satisfies the following
\begin{align*}
\B^{[2]}(e_{12})&=4K_{\text{min}}(x)e_{12}-4\varepsilon\|\a_{13}\|^2e_{34},\\
\B^{[2]}(e_{13})&=\big(K_{12}+K_{34}+2K_{\text{min}}(x)\big)e_{13}-2\varepsilon\|\a_{13}\|^2e_{24},\\
\B^{[2]}(e_{14})&=\big(K_{12}+K_{34}+2K_{\text{min}}(x)\big)e_{14}
+2\varepsilon\|\a_{13}\|^2e_{23},\\
\B^{[2]}(e_{23})&=\big(K_{12}+K_{34}+2K_{\text{min}}(x)\big)e_{23}+2\varepsilon\|\a_{13}\|^2e_{14},\\
\B^{[2]}(e_{24})&=\big(K_{12}+K_{34}+2K_{\text{min}}(x)\big)e_{24}
-2\varepsilon\|\a_{13}\|^2e_{13},\\
\B^{[2]}(e_{34})&=4K_{\text{min}}(x)e_{34}-4\varepsilon\|\a_{13}\|^2e_{12}.
\end{align*}

Then, we easily verify that the 2-vectors $\{\eta_i\}_{1\leq i\leq6}$ 
are eigenvectors of $\B^{[2]}$ with corresponding 
eigenvalues as stated in part $(ii3)$. The 
nonnegativity of the eigenvalues follows 
from Lemma \ref{sb} and part $(i)$ of the proposition.\qed 
\vspace{1ex}

\noindent\emph{Proof of Theorem \ref{thm1} for $n=4$.} 
Part $(i)$ of Proposition \ref{nnic} implies that $M^4$ has 
nonnegative isotropic curvature. We begin by proving 
the theorem for simply connected submanifolds. To this 
aim, we distinguish two cases.
\vspace{0.5ex}

\noindent\emph{Case I.} We assume that there exists a 
point where $M^4$ has positive isotropic curvature. We then claim 
that the manifold $M^4$ is diffeomorphic to $\Sf^4$. Since 
$M^4$ has nonnegative isotropic curvature, it follows from 
Remark (iv) in \cite[p. 623]{seshardi} that $M^4$ 
carries a metric of positive isotropic curvature. 
Given that $M^4$ is simply connected by assumption, 
the main result in \cite{MM} implies that $M^4$ is 
homeomorphic to $\Sf^4$. Furthermore, $M^4$ 
is locally irreducible. If this were not the case, then 
Theorem 3.1 in \cite{MW} would imply that $M^4$ is isometric to 
a Riemannian product of two compact surfaces, which 
leads to a contradiction.

Hence, $M^4$ is locally irreducible, and by Theorem 
2 in \cite{BSacta}, one of the following cases must hold:
\begin{enumerate}[topsep=1pt,itemsep=1pt,partopsep=1ex,parsep=0.5ex,leftmargin=*, label=(\roman*), align=left, labelsep=-0.5em]
\item[(i)] $M^4$ is diffeomorphic to a spherical space form. 
\item[(ii)] $M^4$ is a K\"ahler manifold biholomorphic to 
$\CP^2$.
\item[(iii)] \hspace{0.001ex} $M^4$ is isometric to a compact symmetric space.
\end{enumerate}

Since $M^4$ is homeomorphic to $\Sf^4$, only cases (i) 
and (iii) can occur. Clearly, in case (i), the manifold must be 
diffeomorphic to $\Sf^4$. The simply connected compact 
symmetric four-manifolds that may arise include 
$\Sf^4,\CP^2$, and $\Sf^2\times\Sf^2$. After 
imposing local irreducibility, only $\Sf^4$ and 
$\CP^2$ remain. Thus, in case (iii), the manifold 
$M^4$ must also be diffeomorphic to $\Sf^4$. 
\vspace{0.5ex}

\noindent\emph{Case II.} Now, suppose that there is 
no point where $M^4$ has positive isotropic curvature. 
In this case, Proposition \ref{nnic} implies that equality 
holds in \eqref{*} at every point $x\in M^4$, and there 
exists an oriented orthonormal four-frame 
$\{e^x_1,e^x_2,e^x_3,e^x_4\}\subset T_xM$ as in part $(ii)$ 
of the proposition. Since $M^4$ has nonnegative 
isotropic curvature, Theorem 4.10 in \cite{MW} implies 
that one of the following holds: 
\begin{enumerate}[topsep=1pt,itemsep=1pt,partopsep=1ex,parsep=0.5ex,leftmargin=*, label=(\roman*), align=left, labelsep=-0.5em]
\item[(a)] $M^4$ carries a metric of positive isotropic 
curvature. 
\item[(b)] $M^4$ is diffeomorphic to 
a product $\Sf^2\times\varSigma^2$, where $\varSigma^2$ is 
a compact surface.
\item[(c)] $M^4$ is a K\"ahler manifold biholomorphic 
to $\CP^2$. 
\end{enumerate}

In case (a), since $M^4$ is simply connected, 
it follows from \cite{MM} that $M^4$ is 
homeomorphic to $\Sf^4$. We can then proceed as 
in \emph{Case I} to conclude that $M^4$ is diffeomorphic 
to $\Sf^4$.

We now claim that case (b) cannot occur. Suppose, 
to the contrary, that $M^4$ is as described in (b). Since 
$M^4$ is simply connected, the surface $\varSigma^2$ must be 
diffeomorphic to $\Sf^2$, and therefore $M^4$ is 
diffeomorphic to $\Sf^2\times\Sf^2$. Consequently, 
$\beta_\pm(M)=1$, which implies the existence of a 
nontrivial self-dual harmonic 2-form $\omega_+$ and 
a nontrivial anti-self-dual harmonic 2-form $\omega_-$. 
By the Bochner-Weitzenb\"ock formula, both forms are 
parallel, and $\<\B^{[2]}(\omega_\pm),\omega_\pm\>=0$ 
at every point. From this, it follows that the eigenvalues 
of $\B^{[2]}$, given in part $(ii3)$ of 
Proposition \ref{nnic}, satisfy
$
\mu_1\mu_2\mu_3=\mu_4\mu_5\mu_6=0, 
$
or, equivalently, 
$$
(K_{\text{min}}-\|\a_{13}\|^2)\big(K_{12}+K_{34}
+2\big(K_{\text{min}}+\|\a_{13}\|^2\big)\big)=0,
$$ 
and 
$$
(K_{\text{min}}+\|\a_{13}\|^2)\big(K_{12}+K_{34}
+2\big(K_{\text{min}}-\|\a_{13}\|^2\big)\big)=0
$$
at every point. The first equality either yields 
$K_{\text{min}}=\|\a_{13}\|^2=K_{12}=K_{34}=0$ or 
$K_{\text{min}}=\|\a_{13}\|^2$. In the latter case, the second 
equality implies that $K_{\text{min}}=\|\a_{13}\|^2=0$. 
Together with the pinching equality, \eqref{K12}--\eqref{B2a2e}, 
and the Gauss equation, these equalities imply that
$c=H=0$ and that $f$ is totally geodesic, 
contradicting the compactness of $M^4$. 

Now, we assume case (c) holds. Then, either 
$\beta_+(M)=1$ and $\beta_-(M)=0$, or 
$\beta_+(M)=0$ and $\beta_-(M)=1$. 
We will only treat the former case, as 
the latter one can be handled in a similar 
manner. Thus, there exists a nontrivial 
self-dual harmonic 2-form 
$\omega_+$. The Bochner-Weitzenb\"ock 
formula implies that $\omega_+$ is parallel, 
and $\<\B_+^{[2]}(\omega_+),\omega_+\>=0$ 
at every point. From part $(ii3)$ of 
Proposition \ref{nnic}, it follows that the 
eigenvalues of $\B_+^{[2]}$, given there, 
satisfy
$
\mu_1\mu_2\mu_3=0, 
$
or equivalently
$$
(K_{\text{min}}-\varepsilon\|\a_{13}\|^2)\big(K_{12}+K_{34}
+2\big(K_{\text{min}}+\varepsilon\|\a_{13}\|^2\big)\big)=0
$$
at every point. 

We claim that if $\mu_2=\mu_3=0$ at a point $x$, then 
$c=0$ and $f$ is totally geodesic at $x$.
Indeed, suppose first that $\varepsilon=1$. Since 
$\mu_2=\mu_3=0$, we have 
$
K_{12}=K_{34}=K_{\text{min}}=\|\a_{13}\|^2=0
$ 
at $x$. It then follows from \eqref{K12}, \eqref{B2a1e} and 
\eqref{B2a2e}, and the Gauss equation that $c=0$ 
and that $f$ is totally geodesic at $x$. 
Now suppose that $\varepsilon=-1$. Since \eqref{K12} 
gives $K_{\mathrm{min}}(x)=K_{13}$, the Gauss equation 
together with \eqref{B2a1e} yields
\begin{align*}
K_{12}+K_{34}\,+&\,2\big(K_{\text{min}}(x)-\|\a_{13}\|^2\big)
\\
&=4\big(c+\<\a_{11},a_{44}\>\big)+\|\a_{11}-a_{44}\|^2-4\|\a_{13}\|^2.
\end{align*}
Hence, since $\mu_2=\mu_3=0$ at $x$, the above 
equality together with part $(ii3)$ of Proposition \ref{nnic} implies that 
$$
4\big(c+\<\a_{11},a_{44}\>\big)=-\|\a_{11}-a_{44}\|^2+4\|\a_{13}\|^2. 
$$
Combining this equality with \eqref{B2a3}, we obtain 
$\a_{11}=\a_{44}$ and $\a_{13}=0$. Finally, 
\eqref{B2a1e} and \eqref{B2a2e} imply that $c=0$, and that 
$f$ is totally geodesic at $x$. 

Clearly, there exists a point $x_0$ at which 
$\mu_2(x_0)=\mu_3(x_0)>0$. Indeed, if these eigenvalues 
vanished everywhere, the preceding argument would imply 
that the immersion is totally geodesic in Euclidean space, 
a contradiction. Thus, $\mu_1(x_0)=0$, 
and from part $(ii3)$ of Proposition \ref{nnic}, 
it follows that the kernel of $\B_+^{[2]}$ at $x_0$ 
is spanned by the vector 
$e^{x_0}_1\wedge e^{x_0}_2+e^{x_0}_3\wedge e^{x_0}_4$.

Let $Z$ be the dual to the self-dual harmonic 
form $\omega_+$. Since this form is parallel, 
we may assume, after possibly multiplying 
by a constant such that $\|\omega_+\|=\sqrt{2}$, 
that at the point $x_0$ we have 
$$
Z_{x_0}=e^{x_0}_1\wedge e^{x_0}_2+e^{x_0}_3\wedge e^{x_0}_4. 
$$
Now, consider the almost complex structure 
$J_{x_0}\colon T_{x_0}M\to T_{x_0}M$ given by
$J_{x_0}e^{x_0}_1=e^{x_0}_2$ and $J_{x_0}e^{x_0}_3=e^{x_0}_4$. 
Then, we have
$$
\omega_+(v,w)=\<J_{x_0}v,w\>\quad\text {for all } v,w\in T_{x_0}M. 
$$

Moreover, we define the skew-symmetric endomorphism $J$ 
of the tangent bundle of $M^4$ such that 
$$
\omega_+(X,Y)=\<JX,Y\>\quad\text {for all }X,Y\in\mathcal X(M). 
$$
Clearly, $J$ is parallel because $\omega_+$ is parallel. 
Now we claim that $J$ is orthogonal, i.e., $\|J_xv\|=\|v\|$ 
for every point $x\in M^4$ and any $v\in T_xM$. Indeed, let 
$V$ be a parallel vector field along a curve $c\colon[0,1]\to M$ 
such that $c(0)=x,c(1)=x_0$ and $V(0)=v$. Obviously, 
the vector field $W=JV$ is also parallel along $c$. Using 
that $J_{x_0}$ is orthogonal we have 
$$
\|J_xv\|=\|W(0)\|=\|W(1)\|=\|J_{x_0}V(1)\|=\|V(1)\|=\|V(0)\|=\|v\|,
$$
which proves the claim. Since $J$ is both 
skew-symmetric and orthogonal, it defines an 
almost complex structure that is also parallel. Hence, the 
triple $(M^4,\<\cdot,\cdot\>,J)$ is a K\"ahler manifold. 

In the case $c>0$, we have $\mu_1=0$ and 
$\mu_2=\mu_3>0$ at every point $x$.
Since 
the kernel of $\B_+^{[2]}$ at $x$ is spanned by the vector 
$e^x_1\wedge e^x_2+e^x_3\wedge e^x_4$ and because 
$\|Z_x\|=\sqrt{2}$, it follows that 
$$
Z_x=\pm\left(e^x_1\wedge e^x_2+e^x_3\wedge e^x_4\right).
$$ 
Since $Z$ is the dual to the self-dual 
form $\omega_+$, we have 
$$
\<J_xv,w\>=\<\<Z_x,v\wedge w\>\>\quad\text {for all } v,w\in T_xM. 
$$
Consequently, 
$$
J_xe^x_1=\pm e^x_2 \quad\text {and}\quad J_xe^x_3=\pm e^x_4
$$ 
at every point $x$.

From \eqref{B2a1} and the above equalities, it follows that 
the second fundamental form of the submanifold 
satisfies 
\be\label{f}
\a_f(JX,JY)=\a_f(X,Y)\quad\text{for all } X,Y\in\mathcal X(M).
\ee

In the case $c=0$, the above argument applies to the open 
subset of $M^4$ where $f$ is not totally geodesic. 
By continuity, \eqref{f} holds throughout $M^4$. Such 
immersions are called $(2,0)$-geodesic (see \cite{eft}). 

Assume that $c=0$. It follows from the Codazzi 
equation that the second fundamental form of $f$ 
is parallel (see the proof of Theorem 4 in \cite{f2}). Such 
submanifolds, also called extrinsic symmetric spaces, 
were classified by Ferus in \cite{f2}. The $(2,0)$-geodesic ones 
are the standard embeddings (cf. \cite[p. 82]{f2} or 
\cite[Section 3]{eq}) of the K\"ahler symmetric 
spaces. Since the K\"ahler manifold $M^4$ is 
biholomorphic to $\CP^2$, the immersion $f$ 
is a composition $f=j\circ g$, where $g$ is 
the standard embedding of $\CP_r^2$ of constant 
holomorphic curvature $\frac{4}{3r^2}$ into a 
sphere $\Sf^7(r)$ and $j\colon\Sf^7(r)\to\R^{4+m}$ 
is an umbilical inclusion with $r=\frac{1}{H}$ (see for instance 
\cite[p. 295]{eft}, \cite{L} or \cite{Sk}). 

Now assume that $c>0$. Consider the immersion 
$F=i\circ f$, where $i\colon\Q^{4+m}_c\to\R^{5+m}$ 
is an umbilical inclusion. The second fundamental form 
of $F$ is given by 
$$
\a_F(X,Y)=\a_f(X,Y)-c\<X,Y\>f\quad\text{for all } X,Y\in\mathcal X(M).
$$
Clearly, the immersion $F$ is $(2,0)$-geodesic. As in 
the previous case, $F$ can be written as a composition 
$F=k\circ g$, where $g$ is the standard embedding of 
$\CP_r^2$ of constant holomorphic curvature 
$\frac{4}{3r^2}$ into a sphere $\Sf^7(r)$ and 
$k\colon\Sf^7(r)\to\R^{5+m}$ is an umbilical inclusion 
with $r=\frac{1}{\sqrt{c+H^2}}$. 
Hence, the immersion $f$ is as described in part $(ii)$ 
of the theorem, thereby completing the proof in the case 
where the manifold is simply connected. 
\vspace{0.5ex}

Now, we claim that this is always the case, namely 
that $M^4$ is simply connected. Consider the 
universal covering $\pi\colon\tilde M^4\to M^4$. 
From Lemma \ref{n>4}, we know that $\tilde M^4$ is 
compact, and thus the fundamental group of $M^4$ is finite. 
Moreover, the isometric immersion $\tilde f=f\circ\pi$ 
satisfies \eqref{*}. Hence, we conclude that either $\tilde M^4$ 
is diffeomorphic to $\Sf^4$, or $\tilde M^4$ is isometric 
to the projective plane $\CP_r^2$ of constant 
holomorphic curvature $\frac{4}{3r^2}$ and the immersion 
$\tilde f$ is as described in part $(ii)$ of the 
theorem. 
 
Assume first that $\tilde M^4$ is diffeomorphic to $\Sf^4$. 
Arguing as in \emph{Case I}, we conclude that $M^4$ is locally 
irreducible. Then, by Theorem 2 in \cite{BSacta}, we have 
that $M^4$ is either diffeomorphic to a spherical 
space form, or isometric to a compact symmetric space. 
Clearly, in the former case, the manifold $M^4$ must be 
diffeomorphic to $\Sf^4$, since the only compact oriented 
four-dimensional space form is $\Sf^4$ \cite[Prop. 4.4, p. 166]{dC}. 
In the latter case, $\tilde M^4$ is a compact symmetric 
space diffeomorphic to $\Sf^4$. Hence, $\tilde M^4$ is 
isometric to a round four-sphere, and $M^4$ is a 
spherical space form. As before, $M^4$ must also be 
diffeomorphic to $\Sf^4$. 

Now suppose that $\tilde M^4$ is isometric 
to the projective plane $\CP_r^2$ of constant 
holomorphic curvature $\frac{4}{3r^2}$, with
$r=\frac{1}{\sqrt{c+H^2}}$, and $\tilde f=j\circ g$, where
$g$ is the standard embedding of $\CP_r^2$ into $\Sf^7(r)$
and $j\colon\Sf^7(r)\to\Q_c^{4+m}$ is an umbilical inclusion. 
Hence, $M^4=\hat M^4/\Gamma$, where 
$\Gamma\subset{\rm{Isom}}(\hat M^4)$ is a 
finite group acting freely on $\hat M^4=\CP_r^2$ 
(see Theorem 2.3.16 in \cite{W}). 
The complex projective plane admits no nontrivial 
quotients, even topologically (see \cite{BR-GB} 
or \cite[p.~185]{Besse0}). 
Hence, the covering map $\pi\colon\tilde M^4\to M^4$ 
must be a diffeomorphism, completing the 
proof of the theorem.\qed

\subsection{Examples of submanifolds satisfying condition \eqref{*}}

We now present a method for constructing geometrically 
distinct submanifolds that are diffeomorphic to the sphere 
$\Sf^n$, while also satisfying \eqref{*}. 

\begin{proposition}\label{ell}
Let $f\colon M^n\to\Q_c^{n+1}$, $n\geq4$, $c\geq0$, be a compact 
convex hypersurface (contained in an open hemisphere if $c>0$) 
with principal curvatures $0<\lambda_1\leq\dots\leq\lambda_n$. 
If the global minimum of $\lambda _1$ and the global 
maximum of $\lambda _n$ satisfy 
\be\label{lambda}
(\max\lambda _n)^2 \leq
\begin{cases}
\,\frac{4}{n}\big(c+2(\min\lambda _1)^2\big)&\text{if\, } n\geq 5 \\[1mm]
\,2c+3(\min\lambda _1)^2&\text{if\, } n=4, 
\end{cases}
\ee
then $f$ satisfies condition \eqref{*} at every point. Moreover, 
if the inequality \eqref{lambda} is strict, then so is \eqref{*} 
at every point. 
\end{proposition}
\proof
Given that $H\leq\max\lambda_n$, it follows that 
$$
b(n,H,c)\leq b(n,\max\lambda_n,c).
$$ 
Consequently, the proof follows immediately from 
the inequality 
$$
K_{\text{min}}\geq c+(\min\lambda_1)^2,
$$
which is a direct consequence of the Gauss equation.\qed
\vspace{1.5ex}

A large class of ellipsoids satisfies \eqref{lambda}. 
Consider, for instance, the ellipsoid in $\R^{n+1}$ 
defined by 
$$ 
\frac{x_1^2}{a_1^2}+\cdots+
\frac{x_{n+1}^2}{a_{n+1}^2} = 1,
$$
where $0<a_1\leq\dots\leq a_{n+1}$. A 
straightforward computation shows that the 
minimum and the maximum of the principal 
curvatures of the ellipsoid are $a_1/a_{n+1}^2$ 
and $a_{n+1}/a_1^2$, respectively. It follows 
that condition \eqref{lambda} is satisfied if 
$$
a_{n+1}\leq a_1 \varepsilon(n), \quad4\leq n\leq8, 
$$
where
$$
\varepsilon(n)=
\begin{cases}
\,\left(\frac{8}{n}\right)^{1/6}&\text{if\, } 5\leq n\leq 8 \\[1mm]
\,3^{1/6}&\text{if\, } n=4.
\end{cases}
$$

By the classical Hadamard theorem, or based on 
the result in \cite{dCW} for hypersurfaces in 
a sphere, any convex hypersurface is 
diffeomorphic to $\Sf^n$.

In this way, we construct numerous compact, 
geometrically distinct submanifolds that are 
topologically spheres and strictly satisfy \eqref{*} 
at every point. Furthermore, this strict form is 
preserved under sufficiently small smooth 
deformations of any such example.

\section*{Acknowledgement}
The author wishes to thank the referee for the 
careful and thorough reading of the first draft of 
this paper, which resulted in numerous improvements.

\bigskip

\noindent Theodoros Vlachos\\
University of Ioannina \\
Department of Mathematics\\
45110 Ioannina -- Greece\\
e-mail: tvlachos@uoi.gr

\end{document}